\documentclass{ifacconf}  

\counterwithin*{section}{part}

\usepackage{graphicx}      
\usepackage[comma]{natbib}        

\usepackage{amssymb,dsfont}
\usepackage{amsmath,amsfonts}
\usepackage{tikz}
\usepackage{pgfplots}
\usetikzlibrary{intersections, pgfplots.fillbetween}
\pgfplotsset{compat=1.18} 

\usepackage{enumitem}

\usepackage{cleveref}
\usepackage{caption}
\usepackage{subcaption}
\usepackage{graphicx}

\usepackage{xcolor}

\newtheorem{theorem}{Theorem}[section]
\newtheorem{lemma}{Lemma}[section]
\newtheorem{remark}{Remark}[section]

\newtheorem{ass}{Assumption}
\newtheorem{example}{Example}

\newcommand{\N}{\mathds{N}}
\newcommand{\R}{\mathds{R}}
\newcommand{\Rp}{\R_{\geq0}}

\newcommand{\vp}{\varphi}
\newcommand{\ve}{\varepsilon}

\newcommand{\cC}{\mathcal{C}}

\newcommand{\fM}{f_{\rm max}}
\newcommand{\gM}{g_{\rm max}}
\newcommand{\gm}{g_{\rm min}}

\newcommand{\vpM}{\|\vp\|_\infty}
\newcommand{\vpm}{\inf_{s \ge 0} \vp(s)}

\newcommand{\dd}[2][ ]{\tfrac{\text{\normalfont d}#1}{\text{\normalfont d}#2}}
\renewcommand{\d}{\ \text{d}}

\newcommand{\al}{\left \langle}
\newcommand{\ar}{\right \rangle}

\newcommand{\setdef}[2]{\left\lbrace\ #1\ \left|\ \vphantom{#1} #2\ \right.\right\rbrace}

\DeclareMathOperator*{\esssup}{ess\,sup}

\begin{document}
\begin{frontmatter}
 \title{On derivative-free sample-and-hold control with prescribed performance\thanksref{footnoteinfo}}

\author[author]{Lukas Lanza}
\address[author]{Technische Universität Ilmenau, Institute of Mathematics, Optimization-based Control Group, Ilmenau, Germany
        {\tt\small lukas.lanza@tu-ilmenau.de}}
\thanks[footnoteinfo]{
        This work was supported by Deutsche Forschungsgemeinschaft (DFG, German Research Foundation; Project-ID 471539468),
        and the Carl Zeiss Foundation (VerneDCt -- Project-ID 2011640173)
        }%

\begin{abstract}
A feedback controller is proposed to perform output reference tracking with prescribed performance for nonlinear continuous-time systems of relative degree two. 
The controller is of sampled-data type, i.e., measurements are available only at sampling times~--~a typical situation in real systems when sensors are involved.
Furthermore, only output information is available, i.e., neither the full state nor derivatives can be used for feedback.
A sufficient uniform sampling rate is derived and the control consists of piecewise constant signals on the sampling intervals, i.e., zero-order hold.
Feasibility of the controller and satisfaction of the control objective are rigorously proven.
The controller is illustrated by a numerical example.
\end{abstract}


\end{frontmatter}
\section{Introduction}
We consider the following control problem:
The output~$y$ of an only partially known nonlinear multi-input multi-output continuous-time dynamical system with well-defined relative degree and stable internal dynamics should follow a given reference trajectory~$y_{\rm ref}$ while satisfying the following objectives
\begin{enumerate}[label=(\roman*),ref=(\roman*)]
\item the tracking is performed with prescribed performance, i.e., for all~$t \ge 0$ the tracking error satisfies $\|y(t) - y_{\rm ref}(t)\| < \psi(t)$  for a given boundary function~$\psi$ (feasible~$\psi$ are given in the set~$\Phi$ below), \label{CO_ppc}
\item measurement of the output's derivative is not available, \label{CO_noderiv}
\item measurements of the output~$y$ are available only at discrete time instances, i.e., $y(k \tau)$, $k \in \N$. \label{CO_sampling}
\end{enumerate}
For the sake of clear presentation, in this brief paper we assume the dynamical system to be of relative degree two.

Control objective \ref{CO_ppc}, for systems with relative degree two, is already achieved by the continuously updated \emph{funnel controller}, cf.~\cite{hackl2013funnel,BergLe18a}.
Here, ``continuously updated'' refers to the fact that the controllers use the output signal $y(t)$ for $t \in \Rp$.
For systems with relative degree larger than one the funnel controller not only relies on the continuous availability of the output signal but also of its higher order derivatives.
This aspect, namely the issue of the availability of the output's derivatives was addressed in~\cite{ilchmann2006tracking}, and in~\cite{berger2018funnel,Lanz22}.
In~\cite{ilchmann2006tracking} a backstepping-like procedure is used to design a filter which is then interconnected with the system.
In~\cite{berger2018funnel,Lanz22} a so-called funnel pre-compensator is used to provide auxiliary derivative signals, which are then used by the controller. 
For both approaches it was shown that the conjunction of the filter with the system to be controlled results in an overall system which is amenable for funnel control.
Early results on tracking with prescribed accuracy while only using output information (no derivative information according to~\ref{CO_noderiv}) were achieved in, e.g.,~\cite{miller1991adaptive} for single-input single-output LTI systems.
Similar to funnel control, though different, is \emph{prescribed performance control}, see~\cite{bechlioulis2014low}, which also achieves~\ref{CO_ppc} (tracking with prescribed behaviour).
This controller is designed for systems in strict feedback form and assumes availability of the system state.
To overcome the issue of accessibility to the state and use output information only, in~\cite{dimanidis2020output} a high-gain observer is included in the controller design.
Similarly, in~\cite{chowdhury2019funnel} and~\cite{liu2020adaptive} a high-gain observer is used in a funnel control scheme.
The letter approaches, however, suffer from the fact that the high-gain parameters must be chosen a-priori ``sufficiently large'', in particular, a sufficient lower bound is not provided.

The third aspect, namely that signals are available only at discrete sampling time instances~\ref{CO_sampling}, is related to \emph{redesign of controllers}, cf.~\cite{nevsic2005lyapunov,grune2008continuous,grune2008sampled}.
These techniques aim to design piecewise constant, i.e., Zero-order Hold (ZoH), inputs such that the trajectories of the sampled-data system converges to the continuous-time trajectories if the sampling interval tends to zero. 
In~\cite{lanza2023sampleddata} a ZoH feedback controller is proposed, which achieves objectives~\ref{CO_ppc} and~\ref{CO_sampling}.
This means that the output~$y$ of a continuous-time system with arbitrary relative degree and stable internal dynamics tracks a reference trajectory~$y_{\rm ref}$ with prescribed performance while signals are available only at discrete sampling times.
Unlike the controller redesign, the feedback controller in~\cite{lanza2023sampleddata} is proven to achieve the control objective for a fixed sampling interval of length~$\tau > 0$, i.e., the feasibility analysis does not use the limit~$\tau \to 0$.
However, like the continuously updated funnel controller, the controller in~\cite{lanza2023sampleddata} relies on the availability of higher order derivatives of the system's output at sampling times.

In the present paper we also take the issue of availability of derivatives into account, i.e., we include objective~\ref{CO_noderiv}.
We propose a ZoH feedback controller, similar to that in~\cite{lanza2023sampleddata} which achieves control objectives~\ref{CO_ppc}~--~\ref{CO_sampling}.
This means that the output of an nonlinear continuous-time system with relative degree and stable internal dynamics follows a reference trajectory with prescribed performance while output measurements are available only at discrete sampling times, and no access to the output's derivative is assumed.
Since we do not use high-gain observers to approximate the derivatives and no further filters are involved, this is a novelty in high-gain feedback tracking control with prescribed performance.

To avoid that the main idea is obscured by technicalities, in this brief paper we study the case of systems with relative degree two.
However, following the ideas in~\cite{yuz2005sampled} we are confident that a generalisation of the proposed method to systems with higher relative degree is possible.

The remainder of this paper is structured as follows.
In \Cref{Sec:ProblemStatement} we formulate the control objective, and
introduce the system class under consideration.
In \Cref{Sec:ControllerDesign} we present the reasoning, which eventually leads to  the proposed controller.
\Cref{Sec:MainResult} contains the main result, namely a statement ensuring feasibility of the developed controller.
We demonstrate successful functioning of the derivative-free ZoH feedback controller by means of a numerical example in \Cref{Sec:Example}.
A summering paragraph and an outlook to related future research in \Cref{Sec:ConclusionAndOutlook} concludes this paper.
The proof of the main result is placed in the appendix.

\textbf{Notation.}
$\Rp:=[0,\infty)$, $\al \cdot, \cdot \ar$ is the standard inner product on $\R^n$, $\|x\|:=\sqrt{\al x,x\ar}$ for $x\in\R^n$.
For an interval $I\subset\R$,  $L^\infty(I,\R^n)$ is the space of measurable essentially bounded
functions $f: I\to\R^n$ with norm $\|f\|_\infty=\esssup_{t\in I}\|f(t)\|$. 
$L^\infty_{\text{loc}}(I,\R^n)$ is the space of locally bounded measurable functions.
$W^{k,\infty}(I,\R^n)$ is the Sobolev space of all $k$-times weakly differentiable functions
$f:I\to\R^n$ {with} $f,\dots, f^{(k)}\in L^{\infty}(I,\R^n)$.

\section{Problem statement} \label{Sec:ProblemStatement}
We introduce the class of systems to be investigated, and formulate the control objective.

\subsection{System class}
We consider multi-input multi-output continuous-time nonlinear dynamical systems in input-output form with stable internal dynamics
\begin{align} \label{eq:System}
  &  \ddot y(t) = f(d(t),y(t),\dot y(t),\eta(t)) \!+\! g(d(t),y(t),\dot y(t),\eta(t)) u(t), \nonumber \\
  &  \dot \eta(t) = h(\eta(t),y(t),\dot y(t)),  \eta(0) = \eta^0 \in \R^l,  \\
&     y|_{[-\tau_0,0]} = y^0  \in \cC([-\tau_0,0],\R^m), \nonumber
\end{align}
where~$y$ represents the output of the system and $u$ is the input. 
The value~$\tau_0 > 0$ is the ``memory'' of the system, i.e. the initial condition is given by the initial trajectory~$y^0$, and~$\eta$ is the internal state not seen directly in the output.
The quantity $d \in L^\infty( \Rp ; \R^p)$ is an unknown bounded disturbance, $f \in \cC(\R^p \times \R^m \times \R^m \times \R^l ;  \R^m)$ is a drift term, and the input is distributed by the matrix-valued function $g \in \cC(\R^p \times \R^m \times \R^m \times \R^l ; \R^{m \times m})$. The internal dynamics are governed by $h \in \cC(\R^l \times \R^m \times \R^m; \R^l)$.
\begin{ass} \label{Ass:g_inv}
    The matrix valued function $g \in \cC(\R^p \times \R^m \times \R^m \times \R^l ; \R^{m \times m})$ is strictly positive definite, i.e.,
\[
    \forall x \in \R^{p+2m+l} \ \forall z \in \R^m \setminus \{0\} \, : \  \al z, g({x}) z \ar > 0.
\]
Note that changing the sign in~\eqref{eq:ZoH} accounts for strictly negative definite~$g$.
\end{ass}

\begin{ass} \label{Ass:BIBS}
The internal dynamics are bounded-input bounded-state stable, i.e., for all~$c_0 > 0$ there exists~$c_1>0$ such that
\begin{equation}
\|(\xi,\zeta)\|_\infty \le c_0 \  \Rightarrow  \ \|\eta(\cdot;0,\xi,\zeta)\|_\infty \le c_1,
\end{equation}
where $\eta(\cdot,0,\xi,\zeta)$ denotes the unique maximal solution of the second equation in~\eqref{eq:System}.
\end{ass}
In fact, for $\xi,\zeta \in L^\infty(\Rp,\R^m)$ \Cref{Ass:BIBS} ensures the existence of a unique global solution of the internal dynamics.

\begin{remark}
    Note that for sufficiently smooth functions $\tilde f, \tilde g, \tilde h$, cf.~\cite{schwartz1999global}, system~\eqref{eq:System} is equivalent to a state space system
    \begin{align*}
        \dot x(t) &= \tilde f(x(t),d(t)) + \tilde g(x(t),d(t)) u(t), \ x(0) = x^0 \in \R^{m + l} \\
        y(t) &= \tilde h(x(t)),
    \end{align*}
    via Byrnes-Isidori representation, cf.~\cite{byrnes1991asymptotic}.
    \end{remark}

\subsection{Control objective}
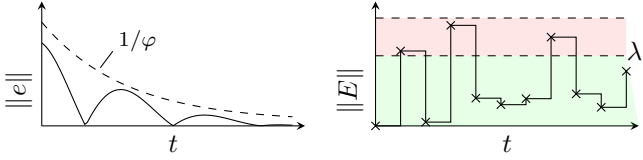
\begin{figure}[ht]
\begin{subfigure}[t]{0.4\linewidth}
\begin{tikzpicture}[>=stealth,scale=0.55]
\draw[->] (0,0) --node[below]{$t$} (6.3,0);
\draw[->] (0,0) -- node[rotate=90,left=1em,above=0.1]{$\|e\|$} (0,3) ;
\draw[scale=2, dashed, domain=0:3, samples=100 , variable=\x, black] plot ({\x}, {1.2*exp(-\x)+0.05});
\draw[scale=2, domain=0:3, samples=100 , variable=\x, black] plot ({\x}, {abs(cos(deg(0.5*\x))*cos(deg(3*\x))*exp(-0.7*\x)});
\draw (1.4,1.4) -- node[above=0.7em,right=0.1em]{\footnotesize{$1/\vp$}} (1.7,1.8);
\end{tikzpicture}
\subcaption{Evolution of the error~$e$ within boundary~$1/\vp$.}
\label{Fig:ErrorInFunnel}
    \end{subfigure}
    \qquad 
    \begin{subfigure}[t]{0.4\linewidth}
    \begin{tikzpicture}[>=stealth,scale=0.55]
\draw[name path= axis,->] (0,0) --node[below]{$t$} (6.3,0);
\draw[->] (0,0) -- node[rotate=90,left=1em, above=0.1em]{$\|E\|$} (0,3) ;
\node (lamb) at (6.2,1.8) {$\lambda$};
\draw[name path=upper, scale=2, dashed, domain=0:3, samples=100 , variable=\x, black] plot ({\x}, 1.3);
\draw[name path=lower, scale=2, dashed, domain=0:3, samples=100 , variable=\x, black] plot ({\x}, 0.65*1.3);
\tikzfillbetween[of=upper and lower]{red, opacity=0.1};
\tikzfillbetween[of=lower and axis]{green, opacity=0.1};
\draw plot[const plot mark left,mark=x,scale=2, domain=0:3, samples=11 , variable=\x, black] ({\x}, {abs(1.5*sin(deg(5*\x))*cos(deg(3*\x))*exp(-0.1*\x)});
\end{tikzpicture}
\subcaption{Activation for~$E(t_k)$ from~\eqref{eq:Ek} with~$\lambda \in (0,1)$.}
\label{Fig:SafeAndSafetycritical}
\end{subfigure}
\caption{Schematic evolution of the tracking error~$e$, and activation profile of~$E(t_k)$.}
\label{Fig:ErrorAndFunnel}
\end{figure}

We aim to design an output feedback ZoH controller, i.e., a feedback law, which, using only sampled output data (no output derivative information), produces a piecewise constant control signal
\begin{equation}
    u(t) = \hat u \quad \forall t \in [t_k, t_k + \tau),
\end{equation}
where~$\tau > 0$ is the length of the sampling interval,
which achieves that the output~$y$ of system~\eqref{eq:System} follows a given reference signal~$y_{\rm ref} \in W^{2,\infty}(\Rp,\R^m)$ with predefined precision in the transient phase.
The latter means that the tracking error $e(t) := y(t) - y_{\rm ref}(t)$ satisfies
\begin{equation} \label{eq:ControlTask}
    \forall \, t \ge 0 \, : \ \vp(t)\| e(t) \| < 1,
\end{equation}
where~$\vp$ belongs to the set
\begin{equation*}
    \Phi:=
         \{\vp\in W^{1,\infty}(\Rp,\R) |
         \inf_{s \ge 0} \vp(s) > 0 \}.
\end{equation*}
The control objective~\eqref{eq:ControlTask} is illustrated in \Cref{Fig:ErrorInFunnel}.

The (possibly time-varying) error tolerance~$\vp \in \Phi$ is chosen by the control engineer, i.e., it is a design parameter.

\section{Controller Design} \label{Sec:ControllerDesign}
We present a ZoH feedback controller, which achieves the control task~\eqref{eq:ControlTask} for systems~\eqref{eq:System} by output feedback only, i.e., the controller does not involve output derivative information.

\subsection{Preliminary calculations}
To propose the control law properly, we record some observations and calculations.
First, for $e=y-y_{\rm ref}$ and $\vp \in \Phi$ we informally define the following auxiliary signals
\begin{equation} \label{eq:ek}
\begin{aligned}
e_1(t) &:=\vp(t) e(t),  \\
e_2 (t) &:= \vp(t) \dot e(t) + \alpha(\|e_1(t)\|^2) e_1(t),
\end{aligned}
\end{equation}
where~$\alpha : [0,1) \to [1, \infty)$ is a bijection chosen by the control engineer;
a suitable choice is $\alpha(s) = 1/(1-s)$.
We emphasise that the signals in~\eqref{eq:ek} are only well-defined on such intervals~$I \subseteq \Rp$,
where $\|e_1(t)\| < 1$ for all~$t \in I$. 
For the sake of concise notation we omit the explicit construction of respective domains but refer to~\cite{lanza2023sampleddata}.
Note that~$e_2$ involves the derivative~$\dot e(t)$, and thus~$\dot y(t)$.
The controller in~\cite{lanza2023sampleddata} directly makes use of~$e_2(t_k)$ at sampling times~$t_k = k \tau$, $k \in \N$.
The idea in the present paper is to substitute the derivative by finite differences and a respective resolvent term.
We aim to find suitable relations between the error's derivative~$\dot e(t)$ and its approximation via finite differences.
From Taylor's theorem using the Lagrange form of the remainder we know for $\tau = t_k - t_{k-1} $
\begin{equation*}
    e(t_{k-1}) = e(t_k) + (t_{k-1}-t_k) \dot e(t_k) + \frac{1}{2} (t_{k-1}-t_k)^2 \ddot e(\sigma),
\end{equation*}
for some~$\sigma \in (t_{k-1},t_k)$.
This yields
\begin{equation*}
\dot e(t_k) = \frac{e(t_k) - e(t_{k-1})}{\tau} + \frac{\tau}{2} \ddot e(\sigma).
\end{equation*}
So we can express~$e_2(t_k)$ in terms of~$e(t_{k-1}), e(t_k)$ and the system dynamics~\eqref{eq:System} by
\begin{equation} \label{eq:e2_k}
\begin{aligned}
    e_2(t_k) &= \vp(t_k) \frac{e(t_k) - e(t_{k-1})}{\tau} + \alpha(\|e_1(t_k)\|^2) e_1(t_k) \\
    & \qquad + \vp(t_k) \frac{\tau}{2} \ddot e(\sigma) \\ 
    &= \vp(t_k) \frac{e(t_k) - e(t_{k-1})}{\tau} + \alpha(\|e_1(t_k)\|^2) e_1(t_k) \\
    &\qquad + \vp(t_k) \frac{\tau}{2} \big( \tilde f(\sigma) + \tilde g(\sigma) u_{k-1} \big),
\end{aligned}
\end{equation}
for some~$\sigma \in (t_{k-1},t_k)$, where with some abuse of notation we use $\tilde f(\sigma) := f(d(\sigma),y(\sigma), \dot y(\sigma),\eta(\sigma))$, and respective for~$\tilde g(\sigma)$.
Here, $u_{k-1}  \in \R^m$ denotes the input vector on the previous interval~$[t_{k-1}, t_{k})$.
Considering the first two addends of~$e_2(t_k)$, we define for $t_{k-1} = t_k - \tau$ the quantity
\begin{equation} \label{eq:Ek}
   \!\!\! E(t_k) := \vp(t_k) \frac{e(t_k) \!-\! e(t_k - \tau)}{\tau} + \alpha(\|e_1(t_k)\|^2) e_1(t_k), 
\end{equation}
which will be used in the ZoH feedback law.
From~\eqref{eq:Ek} it is clear that an initial trajectory on~$[-\tau,0]$ is required; more precise, a measurement at~$t=-\tau$ is needed.

Before we proceed designing the controller, we recall the following two results.
The first lemma gives worst case estimates on the system dynamics when~\eqref{eq:ControlTask} is satisfied, and $\|e_2(t)\| \le 1$.
The precise formulation requires some technical notation.
To avoid these technicalities, we refer to~\cite[Lem.~2.2]{lanza2023sampleddata}, where the statement is given in a general setting.
\begin{lemma} \label{Lem:fMgMgm}
Let~$I \subseteq \Rp$ be an interval.
If {$\|e_1(t)\| < 1$} and $\|e_2(t)\| \le 1$ for all $t \in I$, then there exist constants $\fM,\gM,\gm > 0$ such that
\begin{equation}
\begin{aligned}
\fM & \ge \| f(d,y,\dot y, \eta)|_I \|_\infty, \\
\gM & \ge \| g(d,y,\dot y, \eta)|_I \|_\infty, \\
\gm & \le \frac{\langle \zeta, g(d,y,\dot y,\eta)|_I \zeta \rangle }{\|\zeta\|^2}, \quad \zeta \in \R^m \setminus\{0\},
\end{aligned}
\end{equation}
where $f(d,y,\dot y, \eta)|_I $ means restricting all arguments of~$f$ to~$I$; respectively for~$g$.
\end{lemma}
Note that the existence of the constants is a result, not an assumption.
Later we will make use of these worst case estimates. We emphasise that knowledge of $d,f,g,h$ is not assumed but only worst case estimates.
The second result yields a global upper bound on the auxiliary variable~$e_1$.
For unique $\xi \in (0,1)$ let
\begin{equation} \label{eq:ve1}
\alpha(\xi^2) \xi = 1 + \left\| \frac{ \dot \vp}{\vp} \right\|_\infty ,  \quad 
\ve_1 := \max \{ \|e_1(0)\|, \xi\} < 1.
\end{equation}
Then, we have the following estimate.
\begin{lemma} \label{Lem:e1e2}
Let~$ I \subseteq \Rp$ be an interval, $\|e(0)\|<1$, and~$\ve_1$ as in~\eqref{eq:ve1}.
If for all~$t \in I$ we have $\| e_2(t) \| \le 1$,
then $\|e_1(t)\| \le \ve_1$ for all~$t \in I$.  
\end{lemma}
The proofs of \Cref{Lem:fMgMgm,Lem:e1e2} are provided in~\cite{lanza2023sampleddata}.

\subsection{Design parameters}
To proceed deriving the derivative-free ZoH feedback controller, 
we introduce some control design parameters. 
Invoking \Cref{Lem:fMgMgm,Lem:e1e2} we define the constants
\begin{align} \label{eq:gamKapp0}
  \bar \gamma &:= \big( 2 \alpha'(\ve_1^2) \ve_1^2 + \alpha(\ve_1^2) \big) \big( \alpha(\xi^2) \xi + \alpha(\ve_1^2) \ve_1 \big) , \\
		\kappa_0 &:= \left\| \frac{ \dot \vp}{\vp} \right\|_\infty \! \! \! \!\! ( 1+ \alpha(\ve_1^2) \ve_1) + \|\vp\|_\infty (\fM + \| \ddot y_{\rm ref}\|_\infty) + \bar \gamma. \nonumber
\end{align}
Invoking $1/\vpm = \|\vp\|_\infty$, for 
an activation threshold $\lambda \in (0,1)$ for the quantity~$E(\cdot)$ from~\eqref{eq:Ek}, cf.~\Cref{Fig:SafeAndSafetycritical}.,
we define
\begin{equation} \label{eq:ControlQuantities}
    \begin{aligned}
       \hat \ve &:= \vpM \left( 1 + \alpha(\ve_1^2)\ve_1 \right), \\
        \hat E &:= \vpM \hat \ve + \alpha(\ve_1^2)\ve_1, \\
        \beta & \ge \frac{2 \vpM \kappa_0}{ \gm} \\
        \tilde F &:= \frac{1}{2} \big( \fM + \gM \vpM \beta/\lambda \big), \\
        \kappa_1& := \kappa_0 + \|\vp\|_\infty \beta \gM.
    \end{aligned}
\end{equation}
With the constants given in \eqref{eq:ve1},~\eqref{eq:gamKapp0}, and~\eqref{eq:ControlQuantities} we may now introduce an upper bound on the sampling time $\tau > 0$, namely
\begin{equation} \label{eq:tau}
    \tau \le \min
    \left\{ \begin{array}{l}
     \frac{\vpm \gm \beta - 2 \kappa_0}{\vpM^2 \tilde F   \hat E}, 
      \frac{\kappa_0}{\kappa_1^2} ,  \\
     \frac{1-\lambda}{\vpM (\tilde F + \gM \lambda) + 
  \kappa_0 } 
\end{array}    \right\}.
\end{equation}

\subsection{Derivative-free sampled-data zero-order hold feedback}
We are now in the position to present the ZoH output feedback controller, which applied to a system~\eqref{eq:System} achieves the control objective~\eqref{eq:ControlTask}.
For~$E(\cdot)$ introduced in~\eqref{eq:Ek}, $\beta > 0$ chosen according to~\eqref{eq:ControlQuantities}, $\lambda \in (0,1)$, and a sampling rate not slower than~$1/\tau$ for~$\tau$ satisfying~\eqref{eq:tau}, set
\begin{equation} \label{eq:ZoH}
\!\!  \forall\, t \!\in\! [t_k, t_k + \tau) : u(t) \!=\! \begin{cases}
        -\beta E(t_k), & \!\!\!\| E(t_k) \| < \lambda , \\
        - \beta \frac{E(t_k)}{\|E(t_k)\|^2}, & \!\!\! \|E(t_k)\| \ge \lambda.
    \end{cases}
\end{equation}
Note that in this feedback law only past and current output measurements (no derivatives) are used.
Moreover, \eqref{eq:ZoH} immediately yields an upper bound on the input, namely $\|u\|_\infty \le \beta/\lambda$.

\section{Main Result} \label{Sec:MainResult}
We present the main result of this article, namely feasibility of the control law~\eqref{eq:ZoH}.
To phrase it, \Cref{Thm:ZoHworks} states that the ZoH feedback law~\eqref{eq:ZoH} applied to a nonlinear continuous-time system~\eqref{eq:System} achieves the control objective of output reference tracking with predefined accuracy~\eqref{eq:ControlTask} while the controller receives only sampled output data.
The latter means that no measurements of the output's derivative are required for the controller.
\begin{theorem} \label{Thm:ZoHworks}
Let a reference~$y_{\rm ref} \in W^{2,\infty}(\Rp,\R^m)$ and a funnel function~$\vp \in \Phi$ be given.
Consider a system~\eqref{eq:System} with initial trajectory~$y^0: [-\tau,0] \to \R^m$ such that
\begin{equation} \label{eq:InitialCondition}
y^0(-\tau) = y_{\rm ref}(-\tau), \quad y^0(0) = y_{\rm ref}(0).
\end{equation}
Then the controller~\eqref{eq:ZoH} applied to system~\eqref{eq:System} is initially and recursively feasible, i.e., $\|e_2(t)\| \le 1$ and $\|e_1(t)\| < 1$ for all~$t \ge 0$.
The latter means satisfaction of the control objective~\eqref{eq:ControlTask}.
\end{theorem}
The proof is relegated to the appendix.

Note that the equality initial condition~\eqref{eq:InitialCondition} is artificial to ensure $u_{-1} = 0$.
In applications, the latter can be achieved by choosing the reference appropriately.

\section{Numerical example} \label{Sec:Example}
\begin{figure}[h!]
\begin{center}
\includegraphics[trim=2cm 4cm 5cm 15cm,clip=true,width=4.3cm]{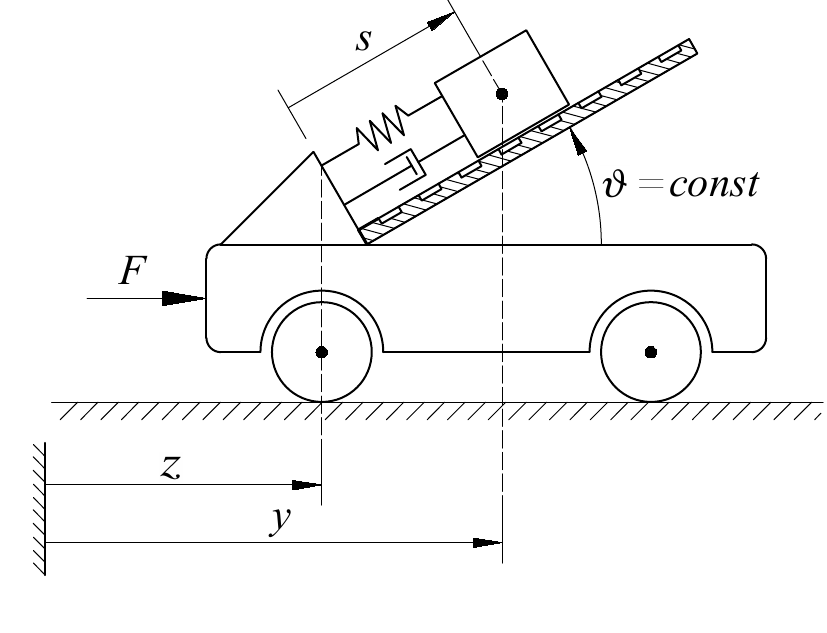}
\end{center}
    \vspace{-5mm}
    \caption{Mass-on-car system. The figure is based on~\cite{SeifBlaj13,BergIlch21}.}
    \label{Fig:MoC}
    \vspace{-2mm}
\end{figure}

We illustrate the proposed controller by a numerical example.
To this end, we consider the widely used mass on car system, cf.~\cite{SeifBlaj13}.
On a car with mass~$m_1$ a second mass~$m_2$ is passively sliding on a ramp.
The second mass is coupled to the car via a spring-damper combination, and the ramp is inclined by an angle~$\vartheta \in (0,\pi/2)$, see \Cref{Fig:MoC}.

The equations of motion in input/output form read
\begin{equation*}
\begin{aligned}
\ddot y(t) &= R_0 y(t) + R_1 \dot y(t) + S \eta(t) + \Gamma u(t), \\
\dot \eta(t) &= Q \eta(t) + P y(t),
\end{aligned}
\end{equation*}
for matrices $R_0,R_1,S,Q,P$ with appropriate dimensions, where~$Q$ is Hurwitz, and $\Gamma > 0$.
We refer to, e.g.,~\cite{BergIlch21} for a detailed derivation of the upper equations.
Note that this system has bounded-input bounded-state stable internal dynamics since~$Q$ is a stable matrix. 
For simulation purpose we choose
$y_{\rm ref} = 0.4 \sin(\pi/2 \, t)$, a constant error tolerance $\psi = 1/\vp = 0.08$, and simulate on the interval~$[0,2]$.
The calculations to obtain the worst case dynamics are performed in~\cite{lanza2023sampleddata} thoroughly, so we use these estimates here.
Starting on the reference we have $e(0)=0$, and we assume that the system satisfies $y(-\tau) = y_{\rm ref}(-\tau)$, i.e., initially $u_{-1} = 0$.

\begin{example} \label{Ex:1}
Choosing~$\lambda=0.7$ conditions~\eqref{eq:ControlQuantities} and~\eqref{eq:tau} are satisfied with~$\beta = 25.2$ and $\tau = 1.8 \cdot 10^{-3}$.
We compare the controller~\eqref{eq:ZoH} with the controller from~\cite{lanza2023sampleddata}, i.e., we compare to the case when the output derivative is available.
Following~\cite[Rem.~5.2]{lanza2023sampleddata}, we use
\begin{equation} \label{eq:ZoH_deriv}
     u_{\rm ZoH + deriv}(t) =  \begin{cases}
        -\beta e_2(t_k), &  \| e_2(t_k) \| < \lambda , \\
        - \beta \frac{e_2(t_k)}{\|e_2(t_k)\|^2}, & \|e_2(t_k)\| \ge \lambda,
    \end{cases}
\end{equation}
where~$e_2$ is given by~\eqref{eq:ek}.
The simulation results are depicted in \Cref{Fig:Signals,Fig:Controls}.
Both controllers~\eqref{eq:ZoH} and~\eqref{eq:ZoH_deriv} achieve the control objective~\eqref{eq:ControlTask} with similar input values.
\begin{figure}
    \centering
    \includegraphics[width=0.75\linewidth]{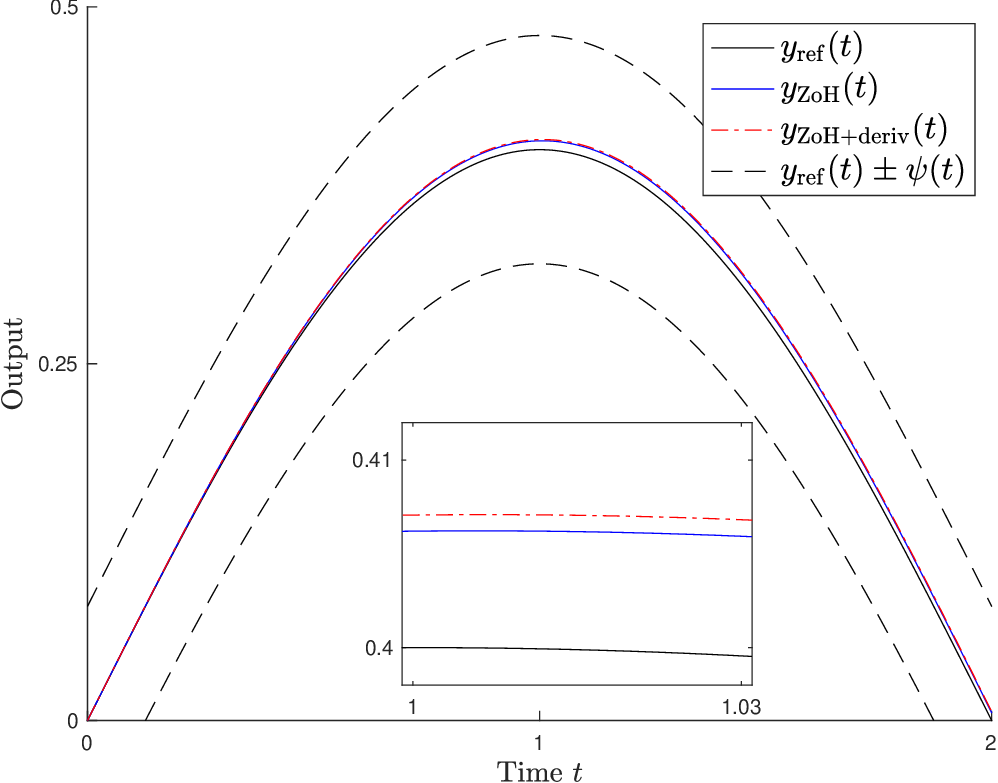}
    \caption{Tracking result for control laws~\eqref{eq:ZoH} and~\eqref{eq:ZoH_deriv}.}
    \label{Fig:Signals}
\end{figure}
\begin{figure}
    \centering
    \includegraphics[width=0.75\linewidth]{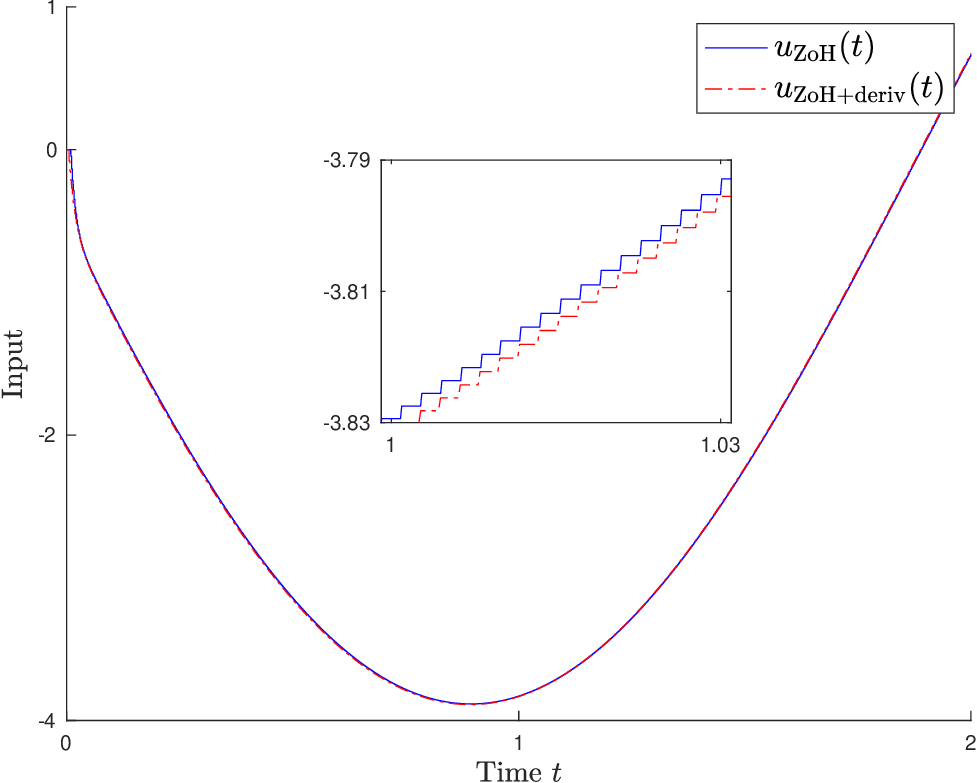}
    \caption{Input signals.}
    \label{Fig:Controls}
\end{figure}
\end{example}

\begin{example}
    The worst case estimates $\fM,\gM, \gm$ from \Cref{Lem:fMgMgm} are very conservative which causes the control parameter~$\beta$ in~\eqref{eq:ControlQuantities} to be large, and the sufficient sampling time~$\tau$ in~\eqref{eq:tau} to be small.
    We illustrate that both controllers from \Cref{Ex:1} also work for less conservative parameters.
    To this end, we run the same simulation as in \Cref{Ex:1} but with~$\beta = 5$ and~$\tau = 0.07$, i.e., the sampling is about forty times slower.
    The simulation results are shown in \Cref{Fig:Signals_LessConservative,Fig:Controls_LessConservative}.
    \begin{figure}
    \centering
    \includegraphics[width=0.75\linewidth]{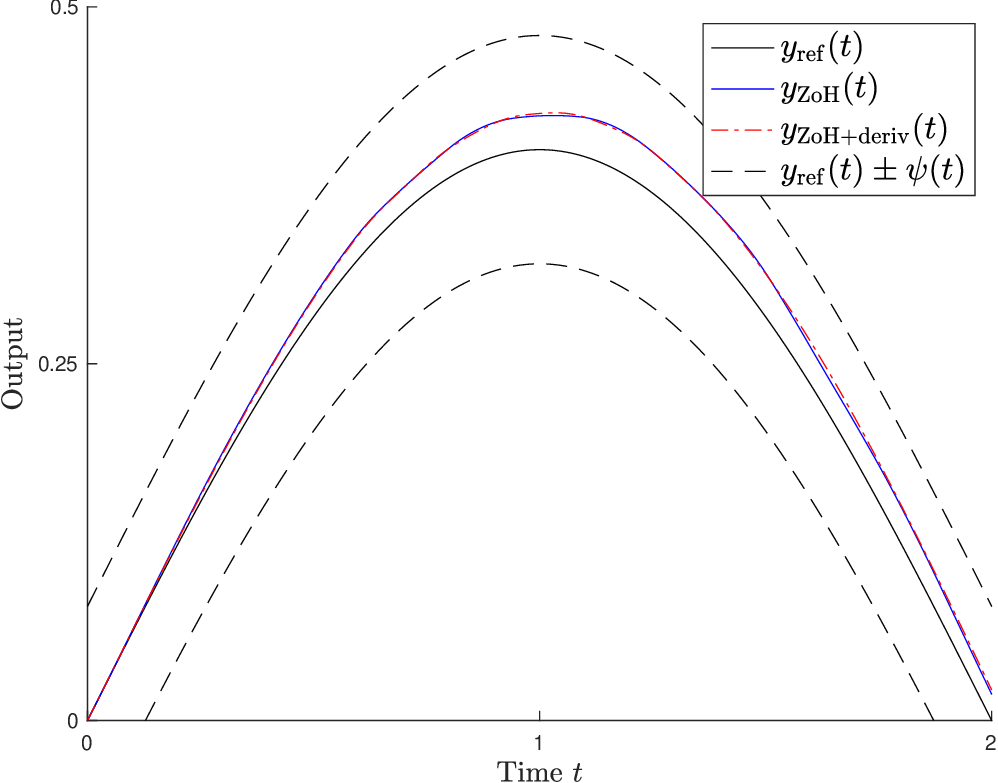}
    \caption{Tracking result for control laws~\eqref{eq:ZoH} and~\eqref{eq:ZoH_deriv}.}
    \label{Fig:Signals_LessConservative}
\end{figure}
\begin{figure}
    \centering
    \includegraphics[width=0.75\linewidth]{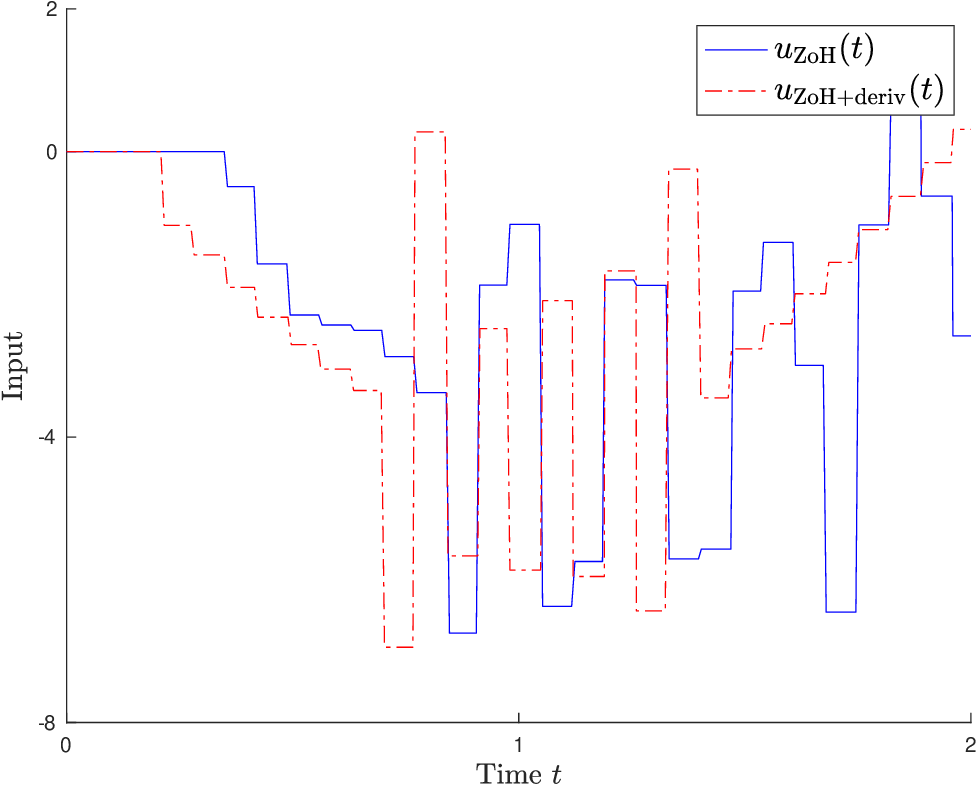}
    \caption{Input signals with slower sampling.}
    \label{Fig:Controls_LessConservative}
\end{figure}
While in \Cref{Ex:1} the control signals for the controller~\eqref{eq:ZoH} and~\eqref{eq:ZoH_deriv} are almost identical, they differ in the second case, see \Cref{Fig:Controls_LessConservative}.
The tracking behaviour in \Cref{Fig:Signals_LessConservative} is still similar for both controllers.
Notably, the tracking performances in \Cref{Fig:Signals} and \Cref{Fig:Signals_LessConservative} are similar.
This further shows that the proposed controller works well even if the sampling frequency is much slower than~$1/\tau$ in~\eqref{eq:tau}.
However, since in this case the conditions of \Cref{Thm:ZoHworks} are not satisfied, no guarantees can be given.
\end{example}

\section{Conclusion and outlook} \label{Sec:ConclusionAndOutlook}
We proposed a sampled-data zero-order hold output feedback controller, which achieves reference tracking with prescribed performance, i.e., the control objective~\eqref{eq:ControlTask} is satisfied.
In contrast to previously proposed controllers,~\cite{hackl2013funnel,BergIlch21,lanza2023sampleddata} the ZoH controller~\eqref{eq:ZoH} does not rely on the availability of output derivative measurement, but uses finite difference approximations.
Feasibility of the proposed controller was proved rigorously.
Since approximation of higher order derivatives is possible with finite differences, we are convinced that a generalisation of the proposed controller to systems with higher relative degree is straightforward.
However, the sampling time~$\tau$ will decrease with increasing relative degree as can be already seen from the analysis in~\cite{lanza2023sampleddata}.
We are confident that finite differences can also be used to 
avoid incorporation of derivative signals in the continuous-time funnel controller~\cite{BergIlch21}.
This will be topic of future research.

\appendix
\section{Proof of Theorem~\ref{Thm:ZoHworks}}
First, we provide some observations and estimates to be used in the main proof.
Invoking \Cref{Lem:e1e2} we may derive the following estimate on the the derivative of the error signal
\begin{align}
\| \dot e \|_\infty &\le \frac{ \sup_{s \ge 0}\| e_2(s)\| + \sup_{s \ge 0} \|\alpha(\|e_1(s)\|^2) \|e_1(s)\|}{\vpm} \nonumber \\ 
& \le \frac{1 + \alpha(\ve_1^2)\ve_1 }{\vpm}   = \hat \ve,
\end{align}
The previous estimate yields an upper bound for the quantity~$E(\cdot)$ in~\eqref{eq:Ek}, independent of~$\tau$, namely
\begin{equation}
\|E\|_\infty \le \vpM \hat \ve + \alpha(\ve_1^2)\ve_1  = \hat E.
\end{equation}
Further, for~$t_k \in \Rp$ we may formally estimate~$e_2(t_k)$ by
\begin{equation} \label{eq:e2_k_estimate}
\| e_2(t_k)\| \le \|E(t_k)\| + \tau  \vpM \tilde F.
\end{equation}

Now we present the proof of \Cref{Thm:ZoHworks}.

\begin{pf}[Proof of \Cref{Thm:ZoHworks}] 
The proof consists of two main steps. 
In the first step we establish the existence of a local solution of the closed-loop initial value problem~\eqref{eq:System},~\eqref{eq:ZoH}.
In the second step we show feasibility of the proposed control law, i.e.,  $\|e_2(t)\| \le 1$, and $\|e_1(t)\| < 1$ for all~$t \ge 0$. 
The latter means that the tracking error evolves within the funnel boundaries.

\noindent
\emph{Step 1.}
The application of the control signal~\eqref{eq:ZoH} to system~\eqref{eq:System} leads to an initial value problem. 
If this problem is considered on the interval~$[0,\tau]$, then there exists a unique maximal solution on $[0,\omega)$ with $\omega\in(0,\tau]$.
If the error variables $e_1(t)$, $e_2(t)$ evolve within the set $\setdef{\xi \in \R^m}{ \| \xi \| < 1}$ for all $t\in [0,\omega)$, then
 $\|(y(\cdot), \dot y(\cdot)) \|$ is bounded on the interval $[0,\omega)$, and by \Cref{Ass:BIBS} $\eta(\cdot,y,\dot y)$ is bounded as well. 
Then $\omega =\tau$, cf. \cite[\S~10, Thm.~XX]{Walt98} and nothing is left to show.
Seeking a contradiction, we assume the existence of $t' \in [0,\omega)$ such that $\| e_\ell(t')\| > 1$ for at least one $\ell \in \{1,2\}$.
Invoking \Cref{Lem:e1e2} it suffices to show 
$\| e_2(t) \| \le 1$ for all $t \in [0,\omega)$.
Before we do so, we record the following observation.
For $\gamma_{1}(t) := \alpha(\|e_{1}(t)\|^2) e_{1}(t)$ we calculate for ${z}(\cdot) := (d(\cdot), y(\cdot), \dot y)(\cdot),\eta(\cdot))$
\begin{equation} \label{eq:J}
\begin{aligned}
    & \dot e_2(t) - \vp(t) g({z}(t)) u = \dot \vp(t) \dot e(t) + \vp(t) \ddot e(t) \\
    & \quad + \dot \gamma_{1}(t) - \vp(t) g({z}(t)) u \\
    &= \frac{\dot \vp(t)}{\vp(t)} (e_2(t) - \gamma_{2}(t)) + \dot \gamma_{1}(t) \\
     &\quad + \vp(t) ( f({z}(t)) - \ddot y_{\rm ref}(t)  ) =: J(t).
\end{aligned}
\end{equation}
\noindent
\emph{Step 2.} We show $\|e_2(t)\| \le 1$ for all $t \in [0,\omega)$.
Since the analysis does not depend on the actual time step~$t_k$, we may set~$t_k=0$,
and investigate the two cases $\| E(0)\| < \lambda$ and $\| E(0)\| \ge \lambda$ (although $\|E(0)\| < \lambda$ by assumption).
 \\
\emph{Step 2.a} We consider $\| E(0)\| < \lambda$. In this case we have $u_0 = - \beta E(0)$,
and possibly $u_{-1} \neq 0$ satisfying $\|u_{-1}\| \le \beta/\lambda$. So~\eqref{eq:e2_k_estimate} is valid for~$t=t_0=0$.
Seeking a contradiction, we suppose the existence of $t^* := \inf \setdef{ t \in (0,\omega)}{ \| e_2(t) \|> 1} $. For the function $J(\cdot)$ introduced in~\eqref{eq:J} we observe $\| J|_{[0,t^*)} \|_\infty \le \kappa_0$ according to \Cref{Lem:e1e2,Lem:fMgMgm}.
Then we calculate for $t \in [0,t^*]$
\begin{align*}
  1 &<  \| e_2(t^*) \| \le \|E(0)\| + \tau \vpM \tilde F  + \textstyle \int_0^{t^*}    \| \dot e_2(s) \|   \d s \\
    & \overset{\eqref{eq:e2_k_estimate}}{\le}  \|E(0)\| + \tau \vpM \tilde F \\
    & \qquad + \textstyle \int_0^{t^*}   \| J(s) \| + \vpM \gM \|E(0)\|  \d s  \\
    &\le  \|E(0)\| + \tau \vpM \tilde F +  \textstyle \int_0^{t^*} \kappa_0 \vpM \gM \|E(0)\|  \d s \\
    & < \lambda + \tau \vpM \tilde F + \kappa_0 \omega + \vpM \gM \lambda \tau < 1,
\end{align*}
where we used $t^* < \omega \le \tau \le (1-\lambda)/(\kappa_0 + \vpM \tilde F + \vpM \gM \lambda )$ by~\eqref{eq:tau}.
This contradicts the definition of $t^*$. \\
\emph{Step 2.b}
We consider $\| E(0)\| \ge \lambda$. 
In this case we have the control $u = - \beta E(0)/\|E(0)\|^2$, and possibly
${u_{-1} \neq 0}$ with $\|u_{-1}\| \le \beta/\lambda$.
Again, we show $\|e_2(t)\| \le 1$ for all $t \in [0,\omega)$.
To this end, seeking a contradiction, we suppose the existence of $t^* = \inf \setdef{(0,\omega)}{ \| e_2(t) \| > 1 }$.
Invoking the initial conditions, in particular $\|e_2(0)\| \le 1$, and continuity of the involved functions, as well as utilising \Cref{Lem:fMgMgm} and~\eqref{eq:J}, we calculate for $t \in [0,t^*]$
    \begin{align*}
       & \dd{t} \tfrac{1}{2} \| e_2(t)\|^2 = \al e_2(t), \dot e_2(t) \ar 
         = \al  e_2(0) + \textstyle \int_0^t  \dot e_2(s) \d s, \dot e_2(t) \ar \\
        & \le \|e_2(0)\| \| J(t) \| + \omega \| \dot e_2|_{[0,t^*]}\|^2_\infty + \vp(t) \al e_2(0), g({z}(t)) u \ar \\
        & \overset{\eqref{eq:e2_k},\eqref{eq:Ek}}{\le}   \| J(t) \|  +  \omega \| \dot e_2|_{[0,t^*]}\|^2_\infty + \tau \vpM^2 \tilde F \beta \|E(0)\|
        \\ & \qquad - \vp(t) \beta \tfrac{\al E(0), g({z}(t)) E(0)\ar}{\|E(0)\|^2} \\
        & \le  \kappa_0 + \omega \| \dot e_2|_{[0,t^*]}\|^2_\infty + \tau \vpM^2 \tilde F \beta \hat E 
        \\ & \qquad - \inf_{s \ge 0} \vp(s) \gm \beta \\
        & \le \kappa_0 + \omega \kappa_1^2 + \tau \vpM^2 \tilde F \beta \hat E  - \inf_{s \ge 0} \vp(s) \gm \beta  \\
        & \le \tau( \kappa_1^2 + \vpM^2 \tilde F \hat E \beta) + \kappa_0 - \vpm \gm \beta  \\
        & < \tau \vpM \tilde F \hat E \beta + 2 \kappa_0 - \vpm \gm \beta < 0
    \end{align*}
  the penultimate line by~$\omega \le \tau$, and the last line by definitions of~$\beta$ and~$\tau$.
  Further, $\| \dot e_2|_{[0,t^*]} \| \le \kappa_1$ was used, and $\| J|_{[0,t^*]} \|_\infty \le \kappa_0$.
    In particular, this yields $\tfrac{1}{2} \dd{t} \| e_2(t)|_{t=0}\|^2 < 0$, by which $t^* > 0$.
    Therefore, we find the contradiction $1 < \| e_2(t^*)\|^2 < \| e_2(0)\|^2 \le 1$.
    Repeated application of the arguments in \emph{Steps 1} and \emph{2} on the interval 
    $[t_k,t_k + \tau]$, $k \in \N$, yields recursive feasibility.
\end{pf}

\begin{ack}
I thank my colleague Dario Dennstädt (U~Paderborn~/ TU~Ilmenau) for fruitful discussions.    
\end{ack}

\bibliography{references}
\end{document}